\documentclass[11pt]{amsart}
\usepackage{amsfonts}
\usepackage{amscd}
\usepackage{amsmath}
\theoremstyle{plain}
\numberwithin{equation}{section}
\newtheorem{theorem}{Theorem}[section]
\newtheorem{corollary}[theorem]{Corollary}

\newtheorem{lemma}[theorem]{Lemma}

\theoremstyle{definition}

\newtheorem{remark}[theorem]{Remark}
\newtheorem{example}[theorem]{Example}

\newcommand{\cA}{\mathcal{A}}

\newcommand{\cF}{\mathcal{F}}

\newcommand{\cM}{\mathcal{M}}
\newcommand{\cO}{\mathcal{O}}
\newcommand{\cH}{\mathcal{H}}

\newcommand{\fa}{\mathfrak{a}}

\newcommand{\fg}{\mathfrak{g}}

\newcommand{\fk}{\mathfrak{k}}

\newcommand{\fn}{\mathfrak{n}}

\newcommand{\fs}{\mathfrak{s}}

\newcommand{\ga}{\alpha}
\newcommand{\gd}{\delta}

\newcommand{\gk}{\kappa}
\newcommand{\gl}{\lambda}

\newcommand{\go}{\omega}

\newcommand{\gt}{\theta}
\newcommand{\gG}{\Gamma}

\newcommand{\gL}{\Lambda}
\newcommand{\gD}{\Delta}
\newcommand{\gO}{\Omega}

\newcommand{\R}{\mathbb{R}}
\newcommand{\C}{\mathbb{C}}

\newcommand{\N}{\mathbb{N}}
\newcommand{\Z}{\mathbb{Z}}

\newcommand{\res}{\mathrm{Res}}
\newcommand{\sign}{\mathrm{sign}}

\newcommand{\Rad}{\mathrm{Rad}}

\newcommand{\vpl}{\varphi_\lambda (m;\cdot )}
\newcommand{\vpa}{\varphi_\lambda (m;a)}
\newcommand{\ACr}{A_\C^{\mathrm{reg}}}
\newcommand{\Ar}{A^{\mathrm{reg}}}
\newcommand{\dA}{d\mu (m; a)}
\newcommand{\dmA}{d\mu (m;\cdot )}
\newcommand{\dn}{d\nu (m;\lambda )}
\newcommand{\dna}{d\nu (m;\cdot )}
\newcommand{\s}{\, :\, }
\newcommand{\Ht}{{\mathcal H}_t(m)}
\newcommand{\AO}{A(\gO )}

\renewcommand{\Re}{\mathrm{Re}}
\renewcommand{\Im}{\mathrm{Im}}

\theoremstyle{plain}

\begin{document}
\title[The Segal-Bargmann transform]{The Segal-Bargmann transform for the heat equation
associated with root systems}
\author{Gestur \'{O}lafsson and Henrik Schlichtkrull}
\address{Department of Mathematics, Louisiana State
University, Baton Rouge,
LA 70803, USA} \email{olafsson@math.lsu.edu}
\address{Matematisk Institut, K{\o}benhavns Universitet,
Universitetsparken 5, DK-2100 K{\o}benhavn {\O}, Denmark}
\email{schlicht@math.ku.dk}
\subjclass{33C67}
\keywords{Heat Equation, Hypergeometric functions, Riemannian symmetric spaces, Hilbert spaces
of holomorphic functions}
\thanks{Research of \'Olafsson was supported by NSF grants DMS-0402068 and
DMS-0139783. 
}
\begin{abstract}
We study the heat equation associated to a multiplicity
function on a root system, where the corresponding
Laplace operator has been defined by Heckman and Opdam.
In particular, we describe the image of the associated
Segal-Bargmann
transform as a space of holomorphic functions.
In the case where the multiplicity
function corresponds to a Riemannian symmetric space $G/K$ of
noncompact type, we obtain a description of the
image of the space of $K$-invariant $L^2$-function on $G/K$
under the
Segal-Bargmann
transform associated to the heat equation on $G/K$,
thus generalizing (and reproving)
a result of B.~Hall and J. ~Mitchell for spaces of complex type.
\end{abstract}
\maketitle

\section*{Introduction}
\noindent
Let $(M,g)$ be a Riemannian manifold.
The heat equation on $M$ is the partial differential equation
$\Delta u(x,t)=\partial_t u(x,t)$, where $\Delta$ is the
Laplace-Beltrami operator. Consider the
corresponding Cauchy problem $\Delta u(x,t)=\partial_tu(x,t)$,
$\lim_{t\searrow 0}u(x,t)=f(x)$, where $f\in L^2(M,dv)$ and where $dv$ is
the volume form. The Cauchy problem is solved by the family of
linear contraction 
maps 
$$H_t: L^2(M,dv)\ni f\mapsto u(\cdot ,t)=e^{t\Delta }f\in 
 L^2(M,dv).$$ Each map $H_t$ 
has the form of an integral operator 
$H_tf(x)=\int_M f(y) k_t(x,y)\, dv(y)$,
where $k_t$ is the \textit{heat kernel}.
The literature contains an extensive study of $H_t$ and its kernel.
We refer to \cite{D90} and \cite{H01} for introduction and an overview.
One of the natural questions is to describe the image of $H_t$.
This can be done either using real analysis or, in several
special cases, complex methods.

Let $M=\R^r$ for a moment, so that we have the Fourier transform and
a natural complexification at our disposal.
It is easy to see that
\begin{equation}\label{eq-real}
H_t(L^2(\R^r))=\{f\in L^2(\R^r)\mid e^{t|\cdot  |^2}\hat{f}\in L^2(\R^r)\}\, .
\end{equation}
But $H_t$ is also a smoothing operator, in fact the image consists of real
analytic functions, which can be extended to entire functions on $\C^r$.
Hence $H_t$ can also be considered as a linear map
$H_t :L^2(\R^r)\to \cO (\C^r)$, sometimes called the 
\textit{Segal-Bargmann transform} (after \cite{B61,S78}).
It should however be noted that originally 
the Hilbert space $L^2(\R^r)$
was replaced by $L^2(\R^r,d\mu_t)$, where $d\mu_t^r(x)=d\mu_t(x)= (4\pi t)^{-r/2}e^{-|x|^2/(4t)}dx$
is the \textit{heat kernel measure}. The reason is, that $\{d\mu_t^r\}$ forms a \textit{projective
family} of probability measures, which then defines a probability measure on the infinite
dimensional space $\R^\infty$.
The image in $\cO (\C^r)$ of the Segal-Bargmann transform is the \textit{Fock space}
\begin{equation}\label{eq-complex}
\cF_t(\C^r):=\{F\in \cO (\C^r)\mid \|F\|_t^2:=(2\pi t)^{-r/2}\int_{\C^r}
|F(x+iy)|^2
e^{-|y|^2/(2t)}\, dxdy<\infty\}
\end{equation}
and $H_t :L^2(\R^r)\to \cF_t(\C^r)$ is a unitary isomorphism, cf.\ \cite{B61,S78}.

The obvious problem in the general case is that there is
no ``natural'' complexification of $M$. 
An important class of spaces, where such
a complexification exists, is the class of Riemannian
symmetric spaces. The first work in this direction was
the article by Hall \cite{H94}.
Here $\R^r$ is replaced by a connected compact semisimple Lie group $U$, and
$\C^r$ is replaced by its complexification $U_\C$. This was put into a
general framework using polarization of a restriction map in \cite{OO96}.
The results of Hall were
extended to compact semisimple symmetric spaces $U/K$
by Stenzel in \cite{S99}. In
this case
the complexification is given by $U_\C/K_\C$. It is important to note that
in the compact case, every eigenfunction of the algebra of invariant differential
operators as well as the heat kernel itself,
extends to a holomorphic function on $U_\C/K_\C$.
This is related to the fact, that each irreducible representation
of $U$ extends to a holomorphic representation of $U_\C$. In the
noncompact case this does not hold, which makes
the situation more complicatied. The natural complexification in this case is
the \textit{Akhiezer-Gindikin domain}  $\Xi\subset G_\C/K_\C$, see
\cite{AG90}. Using results from \cite{ks05} it was  shown
in \cite{kos05} that the image of the Segal-Bargmann transform on
$G/K$ can be identified as a Hilbert-space of holomorphic functions on
$\Xi$. It was also shown that the norm on this space was not
given by a density function. Some special
cases have also been considered in \cite{H04a,H04b},
but without using the Akhiezer-Gindikin domain explicitly. In
particular, in \cite{H04b} B. ~Hall and J. ~Mitchell 
give a description of
the image of the $K$-invariant functions in $L^2(G/K)$ in case
$G$ is complex. The image is in fact isomorphic to
the  Fock space $\cF_t(\fa_\C )$
described in (\ref{eq-complex}), where $\fa$ is the Lie
algebra of a maximal abelian vector subgroup $A$ in~$G$.

Here we shall also study the case of $K$-invariant functions,
but in a generalized setting.
The main  tool in the study of $K$-invariant
functions is the reduction to analysis on $A$ and $\fa$
by restriction. This is made possible by the Cartan decomposition
$G=KAK$. On $A$ the
radial part of the Laplacian is a singular differential
operator with leading part the Laplacian of $A$. The main
ingredients depend on the system of restricted roots and
the multiplicities $m_\ga$, i.e., the dimension of the root spaces.
The setting of analysis on Riemannian symmetric spaces
has been generalized by
Heckman and Opdam by allowing arbitrary values for the multiplicities
$m_\ga$, which are thus replaced by a Weyl-group invariant function
$m : \gD \to \R$, see \cite{hs94}. 
The spherical functions are replaced by
the Heckman-Opdam hypergeometric functions, and the spherical
Fourier transform becomes the hypergeometric Fourier transform. In order
to have a Plancherel decomposition without discrete parts, it is common to assume
that $m_\ga \ge 0$ for all $\ga\in \gD$. 

In this article we
study the heat equation on $A$ corresponding to a non-negative
multiplicity function.  
In particular our results include the description
of the Segal-Bargmann transform on $L^2(G/K)^K$ for a Riemannian
symmetric space of noncompact type. The paper is organized
as follows. In the first section we introduce the necessary
facts from the Heckman-Opdam theory. We discuss two special
cases, first of all the geometric case, where the set up
corresponds to a Riemannian symmetric space of the noncompact
type, and secondly the case where
the root system is reduced and $m_\ga\in 2\N$ for all $\ga$.
The second section starts off by a short discussion of
the Euclidean case $\R^r\simeq \fa$, and then we discuss
the heat equation corresponding to a nonnegative multiplicity
function. In Theorem \ref{th-fv} we give a description of
the image of the Segal-Bargmann
transform similar to (\ref{eq-real}). The main
result, Theorem \ref{th-main}, 
describes the image in terms of holomorphic function
on $\fa_\C$ using the standard Fock space (\ref{eq-complex}).

\section{The Hypergeometric Functions and
Hypergeometric Fourier Transform}
\noindent
In this section we collect basic facts
about the Heckman-Opdam hypergeometric functions,
the corresponding hypergeometric Fourier transform, and
the Plancherel formula. We then specify those results
to our main example of semisimple symmetric spaces and
to the case where all multiplicities are even.
Our standard references are \cite{hs94,Op95}.

\subsection{The Heckman-Opdam hypergeometric functions}
Let $\fa$ be a $r$-dimensional real vector space
with inner product $(\cdot , \cdot )$. Let
$\Delta \subset \fa^*$ be a root system, not necessarily reduced.
We fix a positive system $\gD^+\subset \gD$
and set
\begin{equation}\label{eq-a+}
\fa^+=\{H\in \fa\mid \forall \ga\in\gD\, :\, \ga (H)>0\}\, .
\end{equation}
For $\lambda \in \fa^*$ define
$h_\lambda\in \fa$ by $\ga (H) = (H,h_\lambda)$ for all $H\in \fa$.
We  define an inner product on $\fa^*$ by $(\lambda ,\mu)
=(h_\lambda ,h_\mu)$. We denote the complex
bilinear extension to $\fa_\C^*$ also
by $(\cdot ,\cdot )$. If $\gl\not=0$ let
$H_\gl =\frac{2}{(\gl ,\gl )} h_\gl$.
Then $\ga (H_\ga )=2$. Note that
$H_\gl$ is independent of the inner product $(\cdot ,\cdot )$.
We extend $(\cdot ,\cdot )$ to a $\C$-bilinear
form on $\fa_\C$ and $\fa_\C^*$.

The associated Weyl group, which is
generated by the reflections
$r_\ga (H)= H-\alpha (H)H_\alpha$, is denoted $W$.
The Weyl group acts on $\fa^*$ (and $\fa_\C^*$) by duality:
$w\gl (H)=\gl (w^{-1}H)$.
Recall that $W\fa^+=\{H\in\fa\mid \forall \ga\in \gD\s
\ga (H)\not=0\}$ is open and dense in $\fa$.
Furthermore, if $w\fa^+\cap \fa^+\not=\emptyset$ then
$w=e$, and $\fa\setminus W\fa^+$ has measure zero.

A \textit{multiplicity function} on $\gD$ is a
function $m :\gD\to \C$, invariant under the Weyl group.
We denote by $\cM_+$ the set of non-negative, i.e., $ m(\ga )\ge 0$,
multiplicity functions. If not otherwise stated, then
we will always assume that $m\in \cM_+$. We
write $m_\ga$ for $m(\alpha)$. We note, that our notation differs
from that of Heckman and Opdam as their root system
is $\{2\ga \mid \ga \in \gD\}$. The present notation is adapted
to the special case of a Riemannian symmetric space explained in
Example \ref{e-G/K}, and it is the same as in \cite{OP04}. Let
\begin{equation}\label{eq-AC}
A_\C =\fa_\C/\Z \{i\pi H_\ga\mid \ga\in\gD\}\, .
\end{equation}
Then $A_\C$ is a $r$-dimensional complex torus and $A_\C = AT$, where
$A=\fa$ and $T=i\fa /\Z \{i\pi H_\ga\mid \ga\in\gD\}$ is compact.
Denote by $\exp :\fa_\C\to A_\C$ the cononical projection.

Define
\begin{eqnarray*}
\ACr &=&
\exp \{ H\in \fa_\C\mid \forall \ga\in \gD\, :\, \ga (H)\not\in \pi i\Z\}\, ,\\
\Ar&=&A\cap \ACr=\exp\{H\in\fa\mid \forall \ga\in\gD\s \ga(H)\not=0\}
\end{eqnarray*}
and
$A^+=\exp \fa^+\subset \Ar$.
Then $WA^+=\Ar$ is open and dense in $A$, $A\setminus \Ar$ has
measure zero and $wA^+\cap A^+=\emptyset$ if $w\not= e$.

For $H\in \fa$ denote by $\partial (H)$ the directional derivative
$\partial (H)f(a)=\partial_t f(a\exp (tH))|_{t=0}$.
Let $H_1, \ldots , H_r$ be an orthonormal basis of $\fa$. We
define a Weyl group invariant differential operator $L(m)$ on $\Ar$ by
\begin{equation}\label{Laplace operator}
L(m)=\sum_{j=1}^r\partial (H_j)^2+\sum_{\ga\in\gD^+}
m_\alpha\frac{1+e^{-2\ga}}{1-e^{-2\ga}}
\partial (h_\alpha)\, .
\end{equation}
Note that the first part is just $L_A$, the Laplace operator on $A$.

Let $\Pi=\{\ga_1,\ldots ,\ga_r\}$ be a set of simple roots in $\gD^+$.
Set
$\gG_+:=\{\sum_{j=1}^r n_j\ga_j\mid n_j\in\N_0\}$.
Define $\rho (m)$ as usually by
$2\rho(m)=\sum_{\ga \in \gD^+}m_\ga\, \ga $.
For $\mu \in \gG_+$ such that $(\mu ,\mu -2\gl )\not=0$ for
all $\mu \in \gL\setminus \{0\}$, define $\gG_\mu (m;\lambda)$
inductively by
$$\Gamma_0 (m;\gl)=1$$
and
$$(\mu ,\mu -2\gl )\gG_\mu (m;\lambda )=
2\sum_{\ga\in\gD^+}m_\ga \sum_{\stackrel{k\in\N}{\mu - 2k\ga \in \gG_+}}
\gG_{\mu-2k\ga }(m;\gl )(\mu +\rho (m)-2k\ga-\gl,\ga)\, .$$

For $B\subseteq A$ and $\go \subseteq\fa$ we set $B(\go )=
B\exp i\go$.
Let $\gO:=\{H\in \fa\mid \forall \ga\in\gD\, :\, |\ga (H) |<\pi/2\}$,
then $\exp$ is bijective $\fa+i2\Omega\to A(2\Omega)$.
We set
$
a^\lambda := e^{\lambda (H)}\in \C$
for $a=\exp H\in A(2\Omega)$ and $H\in \fa+i2\Omega$.

The  \textit{Harish-Chandra series} $\Phi_\lambda (m;a)$ is
defined on $A^+(2\gO )$ by
$$\Phi_\gl (m;a)=a^{\gl -\rho(m)}\sum_{\mu\in\gG_+}\gG_\mu (m;\gl )a^{-\mu}$$
for $\lambda\in\fa_\C^*$ such that $\gl_\ga:=\frac12\lambda(H_\alpha)
\notin \Z$
for all $\alpha\in\Delta$.
It follows in the same way as Corollary 4.2.3 in \cite{hs94} that
the sum converges and that $\Phi_\gl (m;\cdot )$ is holomorphic
on $A^+(2\gO)$.
Let $\gD_i^+=\{\ga\in\gD^+\mid 2\ga\not\in\gD\}$.
Define the $c$-function by the Gindikin-Karpelevic formula
\begin{equation}\label{def-cfct}
c(m;\lambda )=\kappa_0\prod_{\ga\in\gD_i^+}\frac{2^{-\gl_\ga}\gG (\gl_\ga )}
{\gG (\frac{1}{2}(\gl_\ga + m_\ga/2+1))\gG(\frac{1}{2}(\gl_\ga +m_\ga /2 +m_{2\ga}))}
\end{equation}
where the constant $\gk_0$ is chosen so that $c (m;\rho (m))=1$.
The \textit{hypergeometric  function} $\varphi_\gl(m;a)$,
associated to the triple $(\fa ,\gD,m)$  is
defined on  $A^+(2\gO)$ by
\begin{equation}
\varphi_\gl (m;a) :=\sum_{w\in W} c(m; w \gl ) \Phi_{w\gl }(m;a)\, .
\end{equation}
Although the Harish-Chandra series is well defined on $A^+(2\gO )$,
the estimates we need can only be found in the literature
for $\AO$. We will therefore
restrict our attention to this domain in the following, though
it might not be maximal as some examples indicate.
\begin{theorem}[Heckman-Opdam]\label{l-est} Let
$m\in \cM_+$. 
Then the following holds:
\begin{enumerate}
\item The function $(\gl,a)\mapsto\vpa$
extends to a holomorphic function on $\fa_\C^*\times\AO$;
\item We have $L (m)\vpl = ((\gl ,\gl )-(\rho (m),\rho (m)))\vpl$
for all $\gl\in\fa_\C^*$.
\item Let $\gl\in\fa_\C^*$. There exists a constant $C>0$
such that for all $a=\exp (H_R+iH_I)\in \AO$  with
$H_R\in\fa$  and $H_I\in\gO$  
\begin{equation}\label{eq-est1}
  |\vpa | \leq C\;
e^{-\min_{w\in W} \Im w\gl(H_I )+\max_{w\in W} w\rho (m) (H_I)+
\max_{w\in W} \Re w\gl(H_R )}
\end{equation}
In particular, for all $a=\exp (H) \in A$ with $H\in\fa$
\begin{equation}\label{eq-est2}
 |\vpa | \leq C \;
e^{\max_{w\in W} \Re w\gl(H )}.
\end{equation}
\end{enumerate}
\end{theorem}
\begin{proof} See \cite{hs94} Corollary 4.2.6 and Corollary 4.3.13.
For the last part see \cite{Op95}, Proposition 6.1.
\end{proof}

\subsection{The Hypergeometric Fourier Transform}
Recall that the action of $W$ on $A$ and $\fa^*$ gives rise to
an action on functions on $A$ and $\fa^*$
by $w\cdot f(a)=f(w^{-1}a)$, respectively $w\cdot F(\gl )=F (w^{-1}\gl )$.
Define a density function
$\gd (m;\cdot )$ on $A$ by
\begin{equation}\label{d-delta}
\gd (m; a):=\prod_{\ga\in\gD^+}|a^\ga -a^{-\ga}|^{m_\ga}\, .
\end{equation}
Then $\gd (m;\cdot )$ is $W$-invariant and
$\gd (m;a)=\prod_{\ga\in\gD^+}( a^\ga -a^{-\ga})^{m_\ga}$ on
$A^+$ and $\dA=\gd (a)da$ is a positive measure on
$A$ (or $A^+$). Denote the corresponding $L^p$-spaces of $W$-invariant functions
$$L^p(A, d\mu_{m})^W\simeq  L^p (A^+,|W| \dA )\, , \qquad f\mapsto f|_{A^+}\,  .$$

Let
$$\fa^*_+=\{\gl\in\fa^*\mid \forall \ga\in\gD^+\, :\, (\gl ,\ga )>0\}
=\{\gl\in\fa^+\mid h_\gl\in\fa^+\}\, .$$
Then $W\fa^*_+$ is open and dense in $\fa^*$, $\fa^*\setminus W\fa^*_+$
has measure zero, and $w\fa^*_+\cap \fa^*_+=\emptyset$ if
$w\not= e$.
On $\fa^*$ we consider
the measure $\dn = |c (m;i \gl )|^{-2}\, d\gl$ and note that
$$L^p(\fa^*, \dn )^W\simeq L^p(\fa^*_+,|W| \dn )\, ,
\qquad f\mapsto f|_{\fa^*_+}\, .$$

For $f\in L^1(A,\dA )^W\cap L^2(A,\dA)^W$
define $\cF (m;f)=\hat{f}:\fa^*\to \C$ by
\begin{equation}\label{d:Fourtr}
\cF (m;f)(\gl ):= \int_A f(a)\varphi_{-i\gl }(a) \dA
=|W| \int_{A^+}f(a)\varphi_{-i\gl }(a) \dA\, .
\end{equation}
Equation (\ref{eq-est2}) in
Theorem \ref{l-est} shows that the integral converges absolutely and
that $\cF(m;f)$ is bounded by $C_1\|f\|_1$.
We call $\cF (m;f)$ for the \textit{hypergeometric Fourier transform
of $f$}
and the linear map $\cF(m)$ the \textit{hypergeometric Fourier transform}.

We have the following important result, cf.\ \cite{Op95},
Theorem 9.13. The Lebesgue measures used on $A$ and $\fa^*$
are assumed to be regularly normalized, see below
(\ref{Eucl Fourier}).

\begin{theorem}\label{th-Fourier} Assume that
$m\in\cM_+$. The operator $\frac1{|W|}\cF$ 
extends to a
unitary isomorphism
$$
L^2(A,  \dmA )^W\simeq
L^2(\fa^*, \dna )^W.
$$
Let  $f\in C_c^\infty (A)^W$, then
\begin{equation}\label{decay of Ff}
\forall N\in\N\, \exists C>0 \, \forall \gl\in\fa^*:
\quad |\cF(m;f)(\gl)| \leq C(1+|\gl|)^{-N}
\end{equation}
and the following Fourier inversion formula holds
\begin{equation}\label{inversion formula}
f(a)= \frac1{|W|^2}\int_{\fa^*}\cF (m;f)(\gl )\varphi_{i\gl }(m;a)
\dn =
\frac1{|W|} \int_{\fa^*_+}\cF (m;f)(\gl )\varphi_{i\gl }(m;a)
\dn.
\end{equation}
Furthermore,
\begin{equation}\label{eq-FL}
\cF (m; L(m)f)(\gl )=-(|\gl |^2+|\rho (m) |^2)
\cF (m; f)(\gl )\,,\quad \gl\in\fa^* .
\end{equation}
\end{theorem}

Note however, that in general
$\cF$ does not map the convolution product
of two functions on $A$ into the product of the Fourier
transforms.

We now relate the hypergeometric Fourier transform to the
usual Fourier transform $\cF_A$ on the abelian group
$A$. Recall that $\cF_A$ is defined on $L^1(A,da)$
by
\begin{equation}\label{Eucl Fourier}
\cF_A(f)(\gl )=\int_A f(a)a^{-i \gl}\, da,
\end{equation}
and that a constant multiple of $\cF_A$ extends to a unitary
isomorphism $L^2(A,da)\simeq L^2(\fa^*,d\gl )$.
The measures are said to be {\it regularly normalized}
when the mentioned constant is 1, as will be assumed here.

For $s,t\in W$ define $c_{s,t}:\fa^* \to \C$ by
$c_{s,t}(m;\gl )= c(m; s^{-1}i\gl )/ c (m;t^{-1}i\gl )$.
Then $|c_{s,t}(m;\gl )|=1$ for all $\gl \in \fa^*$. It
follows that we can define an isometric  map
$\tau_s : L^p(\fa^* ,d\gl )
\to L^p(\fa^* ,d\gl )$ ($s\in W,\, 1\le p\le \infty$) by
\begin{equation}\label{de-tau}
\tau_sF(\gl )=c_{s,e}(m;\gl )F(s^{-1}\gl )\, .
\end{equation}
Note that $\tau_s$ depends on the multiplicity function
$m$ although it is not indicated it in the notation.
A simple calculation show that $\tau_{st}=\tau_s\tau_t$ and hence
we have an action of $W$ on $L^p(\fa^* ,d\gl )$.
As the Fourier transform
is a unitary isomorphism $L^2(\fa^*,d\gl)\simeq L^2(A,da)$, we
can carry this action over to
$L^2(A,da)$. We denote the subspaces of $\tau (W)$-invariant
functions by $L^2(\fa^*,d\gl)^{\tau W}$ and
$L^2(A,da)^{\tau (W)}$.
\begin{lemma}\label{le-intertwining} The multiplication operator
$$L^2(\fa^*,\dn )\ni F\mapsto \Psi_{\fa} (F):=c(m; - 
i\cdot )^{-1} F\in L^2(\fa^* ,d\gl )$$
is a unitary isomorphism such that $\Psi_{\fa} (t\cdot F)=\tau_t\Psi_{\fa}$.
In particular $\Psi_{\fa} :L^2(\fa^*,\dn )^W\to L^2(\fa^* ,d\gl )^{\tau (W)}$
is a unitary isomorphism.
\end{lemma}
\begin{proof} This follows by simple computation.
\end{proof}
Set
\begin{equation}\label{def-Abel}
\cA = \cF_A^{-1}\circ \cF(m): L^2(A,\dA)^W\to L^2(A, da)^W.
\end{equation}
In particular, for $f\in C_c^\infty (A)^W$ it follows from
(\ref{decay of Ff}) that
\begin{equation}\label{Abel of Ccinfty}
\cA f(a)=\int_{\fa^*} \cF(m;f)(\gl)a^{i\gl}\,d\gl.
\end{equation}
Denote  by $\Psi_A$ the pseudo-differential
operator
\begin{equation}\label{LA}
\Psi_A :=\cF_A^{-1}\circ
\Psi_{\fa} 
\circ \cF_A\,
:=\cF_A^{-1}\circ
\frac{1}{c(m; - 
i\cdot )}
\circ \cF_A .
\end{equation}
\begin{lemma}\label{le-1.4} The map
\begin{equation}\label{eq-PhiA}
\gL:=\frac1{|W|}\Psi_A\circ \cA=
\frac1{|W|}\cF_A^{-1}\circ \Psi_{\fa} \circ \cF(m)
\end{equation}
is a unitary  isomorphism $L^2 (A,\dA)^W\to L^2(A,da)^{\tau (W)}$.
Furthermore, if $f\in C_c^\infty (A)^W$, then
\begin{equation}\label{eq-LL}
\gL (L(m)f)=(L_A-|\rho (m) |^2)\gL (f)\, .
\end{equation}
\end{lemma}
\begin{proof} Equation (\ref{eq-PhiA})
is immediate, and the unitarity of $\gL$
follows directly from Lemma \ref{le-intertwining}.
Finally, (\ref{eq-LL})
follows from (\ref{eq-FL}) and the
corresponding Euclidean expression
$\cF_A (L_Af)(\gl )=-|\gl |^2 \cF_A(f) (\gl )$.
\end{proof}

\begin{example}\label{e-G/K}\textbf{(The geometric case).}
The motivating example for the previous notation
and theory is the case where $(\fa ,\gD, m)$ corresponds to
a Riemannian symmetric space of non-compact type.
For that, let $G$ be a connected semisimple
Lie group with finite center and $\theta : G\to G$  a
Cartan involution. Then $K=G^\theta =\{x\in G\mid \theta (x)=x\}$
is a maximal compact subgroup of $G$ and $G/K$ is a
Riemannian symmetric space of the non-compact type. Denote by $\gt$ the
derived involution on $\fg$, the Lie algebra of $G$.
Then
$\fg=\fk\oplus \fs$ where
$\fk=\fg^\theta=\{X\in \fg\mid \gt (X)=X\}$ and
$\fs= \{X\in \fg\mid \gt (X)=-X\}$. Note that $\fs$
can be identified with the tangent space $T_{x_0} (G/K)$,
where $x_0=eK\in G/K$.
Let $\fa$ be a maximal abelian subspace of $\fs$. For $\ga\in\fa^*$
let
$\fg^\ga:=\{X\in\fg\mid \forall H\in \fa\s [H,X]=\ga (H)X\}$
and
$$\gD :=\{\ga\in\fa^*\setminus \{0\}\mid \fg^\ga\not= \{0\}\}\, .$$
Then $\gD$ is a
root system
and
the function $m:\gD \to \R^+$ defined by
$m_\ga := \dim \fg^\ga$ is a positive multiplicity function.
We say that the triple $(\fa,\gD,m)$ is
\textit{geometric} if it corresponds to a Riemannian symmetric space
in this way. In this case, when $G/K$ is fixed,
the multiplicity function $m$ is omitted from the notation.

Let $\gD^+$ be a fixed positive system and put
$\fn:=\bigoplus_{\ga\in\gD^+}\fg^\ga\, .$
Then $\fn$ is a nilpotent Lie algebra.
Let $N=\exp \fn $ and $A=\exp \fa$, then
the multiplication map
$$N\times A\times K\ni (n,a,k)\mapsto nak\in G$$
is an analytic diffeomorphism the inverse of which
is the Iwasawa decomposition
$$G\ni x\mapsto (n (x), a(x), k(x))\in N\times A\times K\, .$$
As before, we set $a^\lambda = e^\lambda (a)=e^{\gl (H)}$ if
$a=\exp H\in A $. Denote by $dk$ the normalized
Haar measure on $K$.
The \textit{spherical functions}
on $G/K$
are given by
$$\varphi_\gl (x)=\int_K a(kx)^{\gl +\rho}\, dk$$
where $\gl\in\fa_\C^*$, with $\varphi_\gl=\varphi_{\gl'}$
if and only if $\gl'\in W\gl$,
see  \cite{HC58a,HC58b} and
\cite{He84}.
The hypergeometric function $\varphi_\gl (m;\cdot )$
associated to $(\fa,\Delta,m)$
is exactly the function  $a\mapsto\varphi_\gl(a)$.

Recall that $G=KAK$, that $KA^+K$ is open, dense
and with complement of measure zero in $G$, and that
$K/M\times A^+\ni (kM,a)\mapsto kax_0$ is a diffeomorphism onto its
image,
see
\cite{He78}, Chapter IX, Theorem 1.1 and Corollary 1.2.
We identify
$A$ with the subset $A\cdot x_0\subset G/K$. It follows
that $K$-biinvariant measurable functions are determined by their restriction
to $A^+$.
We obtain from \cite{He84},
Chapter I, Theorem~5.8:
\begin{lemma}\label{le-LW}
The invariant measure on $G/K$ can be normalized so that
$$\int_G f(x)\, dx=\int_{A }f(a)\gd( a)\, da=|W| \int_{A^+}f(a)\gd(m; a)\, da$$
for all $f\in L^1(G/K)^K$. In particular, the restriction
map $f\mapsto f|_A$ defines a unitary isomorphism
$L^2(G/K)^K\simeq L^2(A, \dA)^W\simeq L^2(A^+, |W| \dA)$.
\end{lemma}
The \textit{spherical Fourier transform} $\cF:L^2(G/K)^K\to L^2(\fa^*,\dn)^W$
is defined by
$$\cF(f)(\gl ):=\int_{G/K}f(x)\varphi_{-i\gl }(x)\, dx\, ,\qquad
f\in L^2(G/K)^K\cap L^1(G/K)^K\, .$$
By Lemma \ref{le-LW} it follows that the spherical Fourier transform is
a special case of the hypergeometric Fourier transform.

For $f\in C_c^\infty (G/K)^K$ the \textit{Abel transform} of
$f$ is defined by
\begin{equation}\label{eq-Abel}
\cA f(a)=a^\rho \int_N f(an)\, dn\, .
\end{equation}
Then $\cA f$ is $W$-invariant,
cf.\ \cite{HC58a}, Lemma 17 or \cite{He84}, Chapter I, Theorem 5.7.
Here we normalize the Haar measure $dn$ on $N$ such
that the measure in Lemma \ref{le-LW} is given
by $dx=a^{-2\rho}dndadk=dadndk$. 
Then
\begin{equation}\label{eq-FA}
\cF=\cF_A\circ \cA,
\end{equation}
see \cite{He84}, p.\ 450 eq.\ (7).
This explains our definition of $\cA$ in (\ref{def-Abel}).
An integral formula as (\ref{eq-Abel})
seems not available in general.

If $L\subset G$ and $f$ is a function defined
on $G$, then we set $\res_L(f)=f|_{L}$.
By \cite{He84}, Chapter II, Corollary 5.11, we have that
$\res_{A}: C^\infty (G/K)^K\to C^\infty (A)^W$
is an isomorphism.
Denote by $\mathbb{D}(G/K)$ the
commutative algebra of $G$-invariant differential operators on $G/K$.
If $D\in \mathbb{D}(G/K)$ then there exists a unique
$W$-invariant differential
operator $\Rad (D)$ on $A$ such that
\begin{equation}\label{eq-Rad}
\res_{A}\circ D  =\Rad (D)\circ \res_{A}\, .
\end{equation}
The differential operator $\Rad (D)$ is called the \textit{radial part of $D$}.
Denote by $L_{G/K}$ the Laplace-Beltrami operator on $G/K$. Then,  according to
\cite{He84}, Chapter II,
Proposition 3.9, we have
\begin{equation}\label{radial part of L}
\Rad_{A} (L_{G/K})=L(m).
\end{equation}
\end{example}

\begin{example} \textbf{(The even multiplicity case).}
We now consider the (not necessarily geometric) special case where
$m_\ga\in 2\Z$ and $2\ga\not\in\gD$ for
all $\ga\in \gD$. We refer to
\cite{OP04} for details. It follows from
the definition
in (\ref{def-cfct})
that
$1/c(m;\gl )$ is a polynomial. In fact we have:
\begin{equation}\label{eq-cfct2}
\frac{1}{c(m;\gl )}= \prod_{\ga\in\gD^+}\prod_{k=0}^{m_\ga /2-1}
\frac{\gl_\ga +k}{\rho (m)_\ga +k}\, .
\end{equation}
We therefore get the following lemma in this case:
\begin{lemma}\label{PsiA}
The operator $\Psi_A$ is a constant coefficient differential operator
on $A$ given by
$$\Psi_A= \kappa (m)\prod_{\ga\in\gD^+}\prod_{k=0}^{m_\ga/2-1}
(- 
\frac12 \partial (H_\ga)+k))\, $$
where the constant $\kappa (m)$ is given by
$$\kappa (m)=\prod_{\ga\in\gD^+}\prod_{k=0}^{m_\ga/2-1}\frac{1}{\rho (m)_\ga +k}\, .$$
\end{lemma}

The following theorem is obtained in \cite{OP04}, Theorem 5.1:
\begin{theorem}[]\label{th-OP} Assume
$m_\ga\in 2\Z$ and $2\ga\not\in\gD$ for
all $\ga\in \gD$.
Let $\psi_\gl :=\sum_{w\in W}e^{w\gl}$. There
exists a $W$-invariant
differential operator $D$ on $A$ with analytic coefficients
such that
$$
\gd (m;a)\varphi_\gl (m;a)=
c(m; \gl )c(m; -\gl ) D \psi_\gl (a)
=\prod_{\ga\in\gD^+}\prod_{k=0}^{m_\ga /2-1}
\frac{k^2-\rho (m)_\ga^2}{k^2- \gl_\ga^2}D \psi_\gl(a) $$
for all $\gl\in\fa_\C^*$. Moreover, these expressions
are holomorphic in $\gl$.
\end{theorem}

Theorem \ref{th-OP}  implies the following inversion
formula for the Abel transform in
terms of $D$ and $\delta (a)$:

\begin{corollary}\label{c-1.16} Assume that
$m_\ga\in 2\N$ and that $f\in C_c^\infty
(A)^W$. Then
$\cA f\in C^\infty(A)$ and
$$D \cA f(a)=|W|\delta (m;a) f(a)\, , \quad a\in A.$$
\end{corollary}
\begin{proof} The smoothness of $\cA f$ follows from (\ref{decay of Ff})
and (\ref{Abel of Ccinfty}). Moreover, we can carry $D$ under the
integral sign in the latter equation. From
$W$-invariance of $\cF(m;f)$ and $D$
we therefore obtain
\begin{eqnarray*}
D\cA f(a)&=&\int_{\fa^*}\cF(m;f)(\gl )D a^{i\gl }\, d\gl\\
&=&\frac1{|W|}\int_{\fa^*}\cF(m;f)(\gl )D\psi_\gl (a)\,d\gl\\
&=&\frac1{|W|}\delta(m;a)
\int_{\fa^*}\cF(m;f)(\gl )\varphi_{i\gl }(m;a)\, \frac{d\gl}{|c(m;i\gl)|^2}\\
&=&|W|\delta(m;a) f(a)
\end{eqnarray*}
where (\ref{inversion formula}) was used in the last step.
\end{proof}

\begin{remark} The differential operator $D$ is
closely related to the
Heckman-Opdam shift operator $D_+(m)$.
We refer to \cite{S03} for discussion on inversion of
the Abel transform using shift operators.
\end{remark}
\end{example}

\begin{example}\textbf{(The geometric case with $G$ complex).}
Let us now consider the simple case where $m_\ga =2$ for all $\ga$.
This corresponds
to the geometric case where $G$ has a complex structure.
Let $\pi (\gl )=\prod_{\ga\in\gD^+} (\ga , \gl)$, then
$c(\gl )=\pi (\rho )/\pi (\gl )$ by (\ref{eq-cfct2}).
It is well
known that 
\begin{equation}\label{eq-varphicomplex}
\varphi_\gl (a)
=c (\gl )\gd (a)^{-1/2}\sum_{w\in W}\sign (w) a^{w\gl}\, .
\end{equation}
where $\gd(a)^{1/2}:=\prod_{\alpha\in\Delta^+}(a^\alpha-a^{-\alpha})$.
Let $\pi (\partial_a)=\prod_{\ga\in\gD^+}\partial (h_\ga)$.
Then
\begin{equation}\label{eq-PsiA}
\Psi_A=
(-1)^{|\Delta_+|}
\pi (\rho )^{-1}\pi (\partial_a)
\end{equation} by
Lemma \ref{PsiA}.

\begin{lemma}\label{formula for D}
Assume that $m_\ga =2$ for all $\ga\in\gD^+$.
Then the differential operator $D$ can be taken as
$$D= \frac{(-1)^{|\gD_+|}}{\pi (\rho )}\gd^{1/2}\pi (\partial_a)$$
Furthermore, with the notation of (\ref{le-intertwining}),
$$D=
\gd^{1/2}\Psi_A\, .$$
\end{lemma}
\begin{proof} We have
$\pi (\partial_a) a^{w\gl}=\sign (w)\pi (\gl )
a^{w\gl}$ and $c(-\gl )=(-1)^{|\gD_+|}c(\gl)$. Hence
\begin{eqnarray*}
(-1)^{|W|}\pi (\rho )^{-1}\gd (a)^{1/2}\pi (\partial_a)\psi_\gl (a)
&=&
(-1)^{|W|}\pi (\rho )^{-1}\gd (a)^{1/2}\pi (\gl )\sum_{w\in W}
\sign (w) a^{w\gl }\\
&=&c (-\gl )^{-1}\gd (a)^{1/2}\sum_{w\in W}\sign (w) a^{w\gl }\\
&=&[c(\gl)c (-\gl )]^{-1}\gd (a)\varphi_\gl (a)\, .
\end{eqnarray*}
The second part follows from (\ref{eq-PsiA}).
\end{proof}
\begin{lemma}\label{W-action in the complex case}
Assume that $m_\ga =2$ for all $\ga$.
The representation $\tau$ is given on
functions on $\fa^*$ and $A$
by
$\tau (w)F(\gl )=\sign (w)F(w^{-1}\gl)$ and
$\tau (w)f(a)=\sign (w)f(w^{-1}a)$ respectively.
\end{lemma}
\begin{proof} This follows from the fact,  that $c(w\gl )=\sign (w)c(\gl )$.
\end{proof}
\begin{lemma} Assume that $m_\ga =2$ for all $\ga$.
The isometry $\gL : L^2(A,d\mu)^W\to L^2(A,da)^{\tau (W)}$
of Lemma \ref{le-1.4} is
given by multiplication with $
\gd^{1/2}
=
\prod_{\ga\in\gD^+}(a^\ga - a^{-\ga}).$
\end{lemma}
\begin{proof}
By the definition of $\gL$ in Lemma \ref{le-1.4}, by
Lemma \ref{formula for D} and by Corollary
\ref{c-1.16},
$$\gL=\frac1{|W|}\Psi_A\circ\cA=
\frac1{|W|}\delta^{-1/2}D\circ\cA
=
\delta^{1/2}.$$
\end{proof}
\end{example}

\section{The Heat Equation}
\noindent
In this section we study the
heat equation on $\Ar$
associated with the operator $L(m)$ of (\ref{Laplace operator}).
We start with the classical case which is well known.
We will identify $A$ and $\fa$ without further comments.

The \textit{heat equation} on $\fa$ is
$
L_\fa u(X,t)=\partial_t u(X,t)
$.
We consider the corresponding Cauchy problem
\begin{equation}\label{eq-heatA}
L_\fa u(X,t)=\partial_t u(X,t)\, ,\qquad \lim_{t\searrow 0}u(X,t)=f(X).
\end{equation}
For $f\in  L^2(\fa,dX)$ (and with the above limit
in $L^2$-sense)
the solution
is given by applying the contraction semigroup
$e^{t\gD}$ to $f$, i.e., $u(X,t)=e^{t\gD}f(X)=h_t*f(X)$
where $h_t (X)=(4\pi t)^{-r/2}e^{-|X|^2/(4t)}$ is the heat kernel.
With the present
normalization of meausures this can  be written as follows
\begin{equation}\label{eq-solutionA}
u(X ,t)=
\int_{\fa^*} e^{-t |\gl |^2}\hat{f}(\gl )e^{i\gl (X)}\, d\gl\, .
\end{equation}

Because of the factor $e^{-t|\gl |^2}$ we can replace $X$ in (\ref{eq-solutionA})
by $X+iY$ and the result is a holomorphic function on $\fa_\C$.
Define a
density function
\begin{equation}\label{def-go}
\go_t^{\fa}(X+iY) =\go_t(X+iY):=(2\pi t)^{-r/2}e^{-|Y|^2/2t}
\end{equation}
and set
\begin{equation}\label{def-HtA}
\cH_t(\fa )=\{F\in \cO (\fa_\C )\mid \int_{\fa_\C}|F(Z)|^2\go_t^\fa (Z)\, dXdY<\infty\}
\, .
\end{equation}
Here $\cO (\fa_\C)$ stands for the space of holomorphic functions on $\fa_\C$.
The following is well known, see \cite{B61,S78}:
\begin{theorem}[Bargmann, Segal]\label{th-BS} Let $t>0$.
The space $\cH_t$ is a Hilbert space with a reproducing kernel.
The solution $u(\cdot,t)$ belongs to $\cH_t(\fa)$ for all $f\in L^2(\fa)$, and
the map
$$H_t^\fa : L^2(\fa )\ni f\to u(\cdot,t)\in \cH_t(\fa)$$
is a unitary isomorphism.
\end{theorem}
The transform
$H^\fa_t:L^2(\fa )\to\cH_t(\fa)$ is called the
\textit{Segal-Bargmann transform}.

\smallskip
Consider now a non-negative
multiplicity function $m$ on a root system
$\Delta$. In analogy we take
$L(m)u(a,t)=\partial_t u(a,t)$
where $a\in\Ar$,
as a definition for the \textit{hypergeometric heat equation}
associated to $m$, and we consider the corresponding Cauchy problem
\begin{equation}\label{heateq}
L(m)u(a,t)=\partial_t u(a,t)\, ,\qquad \lim_{t\searrow 0}u(a,t)=
f(a).
\end{equation}
More precisely, we shall study this problem
with the condition that
$f\in L^2(A,\dA)^W$ where $\dA=\delta(a)da$, and with the
above limit with respect to this $L^2$-space.

\begin{example}
The heat equation on a Riemannian manifold is given by
(\ref{eq-heatA}),
with $L_{\fa}$ replaced by the Laplace-Beltrami operator. In
particular,
if
the Riemannian manifold is a Riemannian symmetric space $G/K$
of the
non-compact type
we obtain
\begin{equation}\label{eq-heateq}
L_{G/K}u(x,t)=\partial_tu(x,t)\, , \quad \lim_{t \searrow 0}u(x,t)=f(x)
\end{equation}
If we assume that  $f$
is $K$-invariant, then obviously the solution
$u(\cdot,t)$ is also $K$-invariant.
Taking the radial part as in Example
\ref{e-G/K} (see (\ref{radial part of L})) we get
exactly the problem (\ref{heateq}) for the associated
multiplicity function $m$, with $f\in L^2(\Ar,\dA)^W$.
\end{example}

The problem (\ref{heateq}) is easily solved by means of the
hypergeometric Fourier transform:
\begin{lemma}\label{le-hsol} Let $f\in L^2(A,\dA )^W$. The solution to the
problem (\ref{heateq}) is given by
\begin{equation}\label{eq-sol}
u(a,t)=\frac1{|W|^2}
\int_{\fa^*}e^{- t(|\gl |^2+|\rho (m)|^2)}
\cF (m;f) (\gl )\varphi_{i\gl}(m;a)
\, \dn\, .
\end{equation}
In particular, the solution $a\mapsto u(a,t)$ is well defined and real analytic on
$A$ and extends to a $W$-invariant holomorphic function on $\AO$.
\end{lemma}
\begin{proof}
Define $u(a,t)$ for $t>0$ by (\ref{eq-sol}). Then
$u(\cdot,t)$ is well defined and real analytic on
$A$ and extends to a $W$-invariant holomorphic function on $\AO$.
This follows from (\ref{eq-est1}), from which a
uniform estimate $|\varphi_{i\gl}(m;a)|\leq Ce^{k|\gl|}$,
$\gl\in\fa^*$, $a\in A(\Omega)$, is easily derived.

It is now easily seen
that $u(\cdot,t)$ solves the
hypergeometric heat equation. Moreover, applying
dominated convergence as $t\to 0$, we see that
the functions
$\gl\mapsto e^{- t(|\gl |^2+|\rho (m)|^2)}
\cF (m;f) (\gl )\varphi_{i\gl}(m;a)$ converge in $L^2(\fa^*,\dna)^W$ to
$\gl\mapsto\cF (m;f) (\gl )\varphi_{i\gl}(m;a)$. Hence
the limit relation in
 (\ref{heateq}) follows by continuity of the inverse
Fourier transform.
\end{proof}

The extension to $A(\Omega)$ of the solution
$u(a,t)$ in (\ref{eq-sol}) is denoted
$H_t(m;f)$, and the map
$$L^2(A,\dA )^W\ni f\mapsto H_t(m;f)\in \cO(\AO)$$
is called the
\textit{Segal-Bargmann transform} associated with
the hypergeometric heat equation.
The following lemma characterizes of the image of the   transform.
\begin{theorem}[The image, first version]\label{th-fv}
The image
$$\Ht:=H_t(L^2(A,\dA )^W)\subset \cO(\AO)$$
is the space of $W$-invariant holomorphic function $F$ on $\AO$
satisfying the following  conditions:
\begin{enumerate}
\item The function $F|_{A}$ is in $L^2(A,\dA)^W$;
\item The function
$\gl \mapsto e^{(|\gl |^2+|\rho (m)|^2)t}\cF(m; F|_A )(\gl )$
is in $L^2(\fa^* ,\dn )^W$.
\end{enumerate}
\end{theorem}
\begin{proof}
Let $f\in L^2(A,da)^W$ and $F=H_t(m;f)$. Then $F$ is
a $W$-invariant holomorphic function on $\AO$,
and by
(\ref{eq-sol}) we have
\begin{equation}\label{F|A}
F|_A=\cF(m)^{-1}(e^{-(|\cdot |^2+|\rho (m)|^2)t}\cF(m;f)).
\end{equation}
Since $|e^{-(|\gl |^2+|\rho (m)|^2)t}|\le 1$
we have $e^{-(|\cdot |^2+|\rho (m)|^2)t}\cF(m;f)\in L^2(\fa^*,\dn)^W$.
Thus $F|_A\in L^2(A,d\mu)^W$ and
$e^{(|\cdot |^2+|\rho (m)|^2)t}\cF(m;F|_A)=\cF(m;f)\in L^2(\fa^*,\dn)^{W},$
so $F$ satisfies (1) and (2).

Now, let $F\in \cO(\AO)$ satisfy (1)-(2), and let
$$g =\cF (m)^{-1}(e^{(|\cdot  |^2+|\rho (m)|^2)t}\cF(m; F|_A ))\in
L^2(A,da)^W.$$
Let $G=H_t(m;g)$. By definition,
$$\cF(m;G|_A)=e^{-t(|\cdot|^2+|\rho(m)|^2)}\cF(m;g)
=\cF(m;F|_A).$$
Hence $G=F$ almost everywhere on $A$.
Since $A$ is a totally
real submanifold of $\AO$
it follows that $G=F$.
\end{proof}

We extend the $\tau$-action of the Weyl group on $L^2(\fa ,d\gl )$
to $\cH_t(\fa)$ by defining it to be
$\tau_w(F):=H^{\fa}_t \tau_w (H^{\fa}_t)^{-1}F$.
If $F=H^{\fa}_tf$, then the action is given
by
\begin{eqnarray*}\tau_sF(X+iY)&=&\int_{\fa }e^{-|\gl |^2t}
c_{s,e}(m;\gl )
\hat{f}(s^{-1}\gl )e^{i\gl (X+iY)}\, d\gl\\
&=&\int_{\fa }e^{-|\gl |^2t}c_{e,s^{-1}}(m;\gl)
\hat{f}(\gl )e^{i\gl (s^{-1}(X+iY))}\, d\gl\, .
\end{eqnarray*}
In particular, the subspace $\cH_t(\fa )^{\tau (W)}$ is defined
by means of this action.

Recall, that we identify $A$ and its Lie algebra $\fa$.
Define a density function on $\fa_\C$ by
\begin{equation}\label{de-density}
\go_t(m;X+iY)=
(2\pi t)^{-r/2}e^{2t|\rho (m)|^2-|Y|^2/2t}\, .
\end{equation}
\begin{lemma}\label{le-extension} Let
$F=H_t (m;f)\in\cH_t(m)$
with $f\in L^2(A ,\dA)^W$. Then
$e^{t|\rho (m)|^2}\gL(F|_A)$ solves  (\ref{eq-heatA})
on $\fa$ with initial value $\gL f\in L^2(\fa ,dX )^{\tau (W)}$.

In particular $\gL(F|_A)$ extends to a $\tau (W)$-invariant
holomorphic function on $\fa_\C$, denoted by $\gL F$, such that
\begin{equation}\label{2-norm}
\|f\|^2=\int_{\fa_\C}
|\gL F (X+iY)|^2 \go_t(m;X+iY)\, dXdY\, .
\end{equation}
\end{lemma}
\begin{proof}
The following chain of
equalities is easily obtained from the
definitions of $\gL$, $H_t$ and $H^{\fa}_t$:
\begin{eqnarray*}
\cF_A(e^{t|\rho (m)|^2} \gL(F|_A))
&=&
e^{t|\rho (m)|^2} {c(m; -
i\cdot)}^{-1}\cF(m;F|_A)\\
&=&
e^{-t|\gl|^2}{c(m; - 
i\cdot)}^{-1}\cF(m;f)\\
&=&
e^{-t|\gl|^2}\cF_A(\gL f)\\
&=&
\cF_A(H^{\fa}_t\gL f).
\end{eqnarray*}
Hence
$
e^{t|\rho (m)|^2} \gL(F|_A)
=H^{\fa}_t\gL f,
$
which is exactly the first statement of the lemma.

The $\tau(W)$-invariance of $\gL F$ is immediate from
the above definition, since $\gL f$ is $\tau(W)$ invariant by
Lemma \ref{le-1.4}. Finally, (\ref{2-norm}) follows from
the unitarity in Theorem \ref{th-BS}.
\end{proof}

We can now define
a sesquilinear product on $\Ht$ by
\begin{equation}
(F,G):=\int_{\fa_\C}\gL F(X+iY)\overline{\gL G(X+iY)}\,
\go_t(m;X+iY)\, dXdY\, .
\end{equation}
We collect the main results in the following theorem.
\begin{theorem}[The image, second version]\label{th-main} The
space $\Ht$  is a Hilbert space
and $H_t(m):L^2(A,da)^W\to \Ht$ is a unitary isomorphism.
\end{theorem}
\begin{proof}
Follows immediately from the preceding lemmas.
\end{proof}

\begin{example}\textbf{(The geometric case with $G$ complex).}
Assume $m_\ga=2$ for all $\ga$, i.e., $(\fa,\gD,m)$ corresponds
to a Riemannian symmetric space $G/K$ with $G$ complex.
Note that the function $\delta^{1/2}$ has a holomorphic extension
to $A$. We then obtain the following
result of B. Hall and J. Mitchell, \cite{H04b} Theorem 3.

\begin{theorem} Let $f\in L^2(G/K)^K$, and let
$u(x,t)=H_tf(x)$ be the solution to (\ref{eq-heateq}).
The map $X\mapsto u(\exp X,t)$, $X\in\fa$ has a meromorphic
extension to $\fa_\C$, denoted $Z\mapsto U(Z)$,
which is $W$-invariant and satisfies the following:

The map $\delta^{1/2}U$ is holomorphic and
$$\|f\|^2
=(2\pi t)^{-r/2}
\int_{\fa_\C} |(\delta^{1/2}U)(X+iY)|^2 e^{2t|\rho|^2-|Y|^2/2t}\,dX\,dY.
$$
Conversely, any meromorphic function $U(Z)$ which is invariant
under $W$ and which satisfies
$$
\int_{\fa_\C} |(\delta^{1/2}U)(X+iY)|^2 e^{2t|\rho|^2-|Y|^2/2t}\,dX\,dY
<\infty
$$
is the Segal-Bargmann tranform $H_tf$ for some $f\in  L^2(G/K)^K$.
\end{theorem}

\begin{proof} Immediate from Theorem \ref{th-main}
by Corollary
\ref{c-1.16}. The statements about $W$-invariance follows from
Lemma \ref{W-action in the complex case},
since $\delta^{1/2}U$ is invariant under the action mentioned there
if and only if $U$ is invariant for the ordinary action.
\end{proof}
\end{example}

\end{document}